# Complete Pascal Interpolation Scheme For Approximating The Geometry Of A Quadrilateral Element*

Sulaiman Y. Abo Diab,

*Department of Structural Mechanics, Tishreen University, Faculty of Civil Engineering, Lattakia, Syria*

e-mail: sabodiab@tishreen.edu.sy

**Abstract**   This paper applies a complete parametric set for approximating the geometry of a quadrilateral element. The approximation basis used is a complete Pascal polynomial of second order with six free parameters. The interpolation procedure is a natural interpolation scheme. The six free parameters are determined using the natural coordinates of the four nodal points (vertices) of the quadrilateral element and the two intersections points of the lines crossing every two opposite edges (poles). The presented scheme recovers the well known Lagrangian interpolation scheme, when every two opposite edges are parallel.  A third order Pascal interpolation scheme is also presented. The four midpoints of the four edges in addition to the six nodal point from the second order case are used as significant nodal points. It is expected to reflect the geometry properties better since the shape functions are complete.

## 1   INTRODUCTION

The lack of convergence of the quadrilateral elements used by finite element method compared with rectangular and parallelogram elements  was observed in an earlier time point [1]. The Lagrangian interpolation scheme is a widely used interpolation scheme in the scope of numerical and computational methods of applied engineering sciences. Lot of text books and research papers, see for instance [2,...,12] , have used this scheme for approximating the differential geometry properties of a quadrilateral shape and establishing a transformation relation between natural and Cartesian coordinates. Lot of efforts have been made to improve the quality of the approximation through applying adaptive procedures, studying mesh distortion as well as

*This paper is dedicated to Professor  Udo. F. Meißner, his support between 1990-1992 is gratefully acknowledged

analysing approximations error. Example reference studies can be found in [13, ...,18]. A survey about construction of shape functions on convex polygons can be found in [19]. About the origin of the finite element method and the application of the isoparametric transformation, the reader is referred to [20],[21].

A complete Pascal interpolation scheme for approximating the geometry of quadrilateral element is applied in the present paper. The presented scheme includes all the terms of the Pascal triangle of certain order. It recovers the Lagrangian scheme in the special case, when every two opposite edges are parallel. The Lagrangian scheme seems to be unable to approximate the general quadrilateral shape or it can approximate it only incompletely. More information than usually used in the Lagrangian scheme must be included in the interpolation procedure. The basic idea of the present application for a quadrilateral element with straight edges is to find six significant nodal points for the use in the interpolation such that, a second order Pascal polynomial with six undetermined parameters can be chosen. For this purpose, the four element vertices (intersection point between every two neighboured edges) and the two poles, where every two opposite edges meet, are selected. The two poles belong to the geometry of the quadrilateral element as well. They play a key role in the approximation. The undetermined parameters are given their geometrical meaning as usual by substituting the special natural coordinates of the significant nodal points and inverting the resulting linear system of equations. It can be shown, that the inversion is possible and leads to shape functions, which include all the second order terms. The present Pascal interpolation scheme recovers the Lagrangian interpolation scheme in the special case, when every two opposite edges are parallel. Therefore, the presented scheme can be adopted directly in a finite element program for approximating the geometry of rectangular and parallelogram shapes. For a 2D third order Pascal interpolation scheme, the four midpoints of the edges in addition to the six nodal point from the second order case are used as significant nodal points. For this case a complete ten term Pascal polynomial can be chosen.

## 2  2D SECOND ORDER INTERPOLATION SCHEME

Whenever is not pointed out, Latin indices range over the Cartesian co-ordinates and indices between round brackets identify the nodal points. For example $i$ ranges over $x^i$ $(i=1,2)$, where $(p)$ denotes the number of the nodal points.

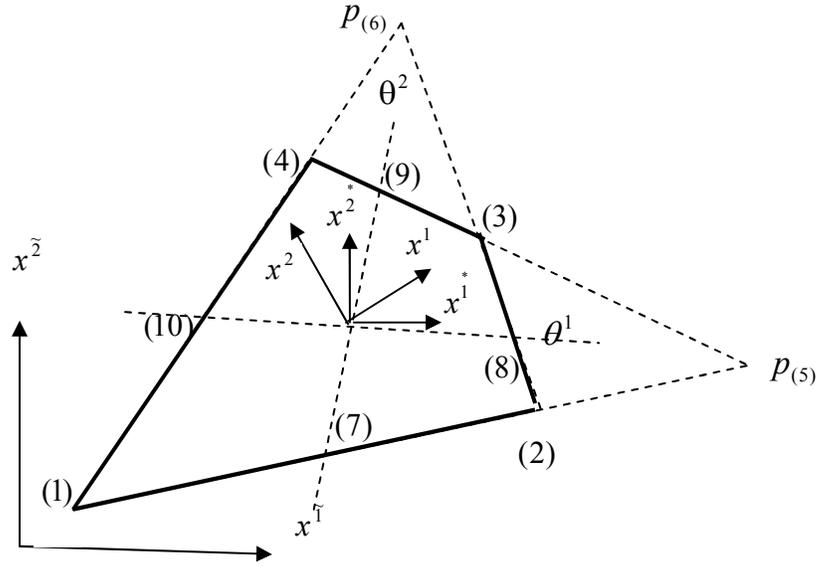

Fig.1 Quadrilateral element, Cartesian and natural coordinate systems

Let $\Omega$ be a quadrilateral domain related to a Cartesian coordinate system $(x^{\tilde{1}},x^{\tilde{2}})$ with the unit vectors $(e_{x^{\tilde{1}}},e_{x^{\tilde{2}}})$ and defined by its four vertices $^{(1),\,(2)\,(3),\,(4)}$ (nodal points)

$$x^{\tilde{i}}_{(p)} = \begin{bmatrix} x^{\tilde{1}}_{(1)} & x^{\tilde{2}}_{(1)} \\ x^{\tilde{1}}_{(2)} & x^{\tilde{2}}_{(2)} \\ x^{\tilde{1}}_{(3)} & x^{\tilde{2}}_{(3)} \\ x^{\tilde{1}}_{(4)} & x^{\tilde{2}}_{(4)} \end{bmatrix} \quad (1)$$

Let the following values define the same nodal points in a natural coordinate system $(\theta^1,\theta^2)$

$$\theta^{i}_{(p)} = \begin{bmatrix} \theta^{1}_{(1)} & \theta^{2}_{(1)} \\ \theta^{1}_{(2)} & \theta^{2}_{(2)} \\ \theta^{1}_{(3)} & \theta^{2}_{(3)} \\ \theta^{1}_{(4)} & \theta^{2}_{(4)} \end{bmatrix} = \begin{bmatrix} -1 & -1 \\ +1 & -1 \\ +1 & +1 \\ -1 & +1 \end{bmatrix} \quad (2)$$

The well-known Lagrangian interpolation approach involves polynomial with four undetermined parameters. Then, the scaled values eqns. (2) of the natural coordinates at the four nodal points are sufficient to establish a transformation relation between the Cartesian coordinates and the natural coordinates.

The current approach starts with a trial Pascal polynomial of second order with six undetermined parameters for approximating the relation between the Cartesian coordinates and the natural coordinates:

$$x^{\tilde{i}} = M^{(p)} a^{\tilde{i}}_{(p)} \qquad (3)$$

In eq. (3), $M^{(p)}$ is a matrix of known functions of the natural co-ordinates, $a^{\tilde{i}}_{(p)}$ are here the total number of undetermined parameters, $x^{\tilde{i}}$ $(i = 1,2)$ are the Cartesian coordinates to be approximated. The detailed form of eq.(3) is given by:

$$\begin{bmatrix} x^{\tilde{1}} & x^{\tilde{2}} \end{bmatrix} = \begin{bmatrix} 1 & \theta^1 & \theta^2 & (\theta^1)^2 & \theta^1\theta^2 & (\theta^2)^2 \end{bmatrix} \begin{bmatrix} a^{\tilde{1}}_{(1)} & a^{\tilde{2}}_{(1)} \\ a^{\tilde{1}}_{(2)} & a^{\tilde{2}}_{(2)} \\ a^{\tilde{1}}_{(3)} & a^{\tilde{2}}_{(3)} \\ a^{\tilde{1}}_{(4)} & a^{\tilde{2}}_{(4)} \\ a^{\tilde{1}}_{(5)} & a^{\tilde{2}}_{(5)} \\ a^{\tilde{1}}_{(6)} & a^{\tilde{2}}_{(6)} \end{bmatrix} \qquad (4)$$

This six parametric approximation functions relate an arbitrary point with the natural coordinates $p'(\theta^1, \theta^2)$ to a corresponding point with the coordinates p $(x^{\tilde{1}}, x^{\tilde{2}})$.

For determining the six free parameters depending on the scaled natural values, not only the values given in eq. (3) are necessary but also the scaled values of additional two other significant points. The two poles: $p_{(5)}$, where the extended lines coincide with the two opposite edges [1] [2] and [3] [4] meet, as well $p_{(6)}$, where the extended lines coincide with the other two opposite edges [2][3] and [4], [1] meet, are chosen as significant, see fig.1. The Cartesian coordinates of poles $p_{(5)}$ and $p_{(6)}$ can be computed by solving the system of the two equations representing the two edges, simultaneously, The equations of the edges [1] [2] and [3] [4] are as follows [22]

$$\begin{vmatrix} x^{\tilde{1}} & x^{\tilde{2}} & 1 \\ x^{\tilde{1}}_{(1)} & x^{\tilde{2}}_{(1)} & 1 \\ x^{\tilde{1}}_{(2)} & x^{\tilde{2}}_{(2)} & 1 \end{vmatrix} = 0 \qquad (5)$$

$$\begin{vmatrix} x^{\tilde{1}} & x^{\tilde{2}} & 1 \\ x^{\tilde{1}}_{(3)} & x^{\tilde{2}}_{(3)} & 1 \\ x^{\tilde{1}}_{(4)} & x^{\tilde{2}}_{(4)} & 1 \end{vmatrix} = 0 \qquad (6)$$

And for the other two opposite edges $^{(2)(3)}$ and $^{(4),\,(1)}$ we have

$$\begin{vmatrix} x^{\tilde{1}} & x^{\tilde{2}} & 1 \\ x^{\tilde{1}}_{(2)} & x^{\tilde{2}}_{(2)} & 1 \\ x^{\tilde{1}}_{(3)} & x^{\tilde{2}}_{(3)} & 1 \end{vmatrix} = 0 \qquad (7)$$

$$\begin{vmatrix} x^{\tilde{1}} & x^{\tilde{2}} & 1 \\ x^{\tilde{1}}_{(4)} & x^{\tilde{2}}_{(4)} & 1 \\ x^{\tilde{1}}_{(1)} & x^{\tilde{2}}_{(1)} & 1 \end{vmatrix} = 0 \qquad (8)$$

As a result, we get the Cartesian coordinates $(x^{\tilde{1}}_{(5)}, x^{\tilde{2}}_{(5)}), (x^{\tilde{1}}_{(6)}, x^{\tilde{2}}_{(6)})$ of $p_{(5)}$ as well $p_{(6)}$. Now, extending the matrix $x^{\tilde{i}}_{(p)}$ to include these values gives :

$$x^{\tilde{i}}_{(p)} = \begin{bmatrix} x^{\tilde{1}}_{(1)} & x^{\tilde{2}}_{(1)} \\ x^{\tilde{1}}_{(2)} & x^{\tilde{2}}_{(2)} \\ x^{\tilde{1}}_{(3)} & x^{\tilde{2}}_{(3)} \\ x^{\tilde{1}}_{(4)} & x^{\tilde{2}}_{(4)} \\ x^{\tilde{1}}_{(5)} & x^{\tilde{2}}_{(5)} \\ x^{\tilde{1}}_{(6)} & x^{\tilde{2}}_{(6)} \end{bmatrix} \qquad (9)$$

The matrix $\theta^{i}_{(p)}$ must be also extended to include the scaled values of the natural coordinates of the points $p_{(5)}$ and $p_{(6)}$. For this purpose, the position vectors of the midpoints of the two element edges $^{(2)(3)}$ and $^{(4),\,(1)}$ are calculated as follows:

$$\vec{r}_{(8)} = (\frac{x^{\tilde{i}}_{(2)} + x^{\tilde{i}}_{(3)}}{2})\vec{e}_{\tilde{i}}; \qquad (10)$$

$$\vec{r}_{(10)} = (\frac{x^{\tilde{i}}_{(4)} + x^{\tilde{i}}_{(1)}}{2})\vec{e}_{\tilde{i}}$$

and the for other two edges $^{(1)\,(2)}$, $^{(3)\,(4)}$ as follows

$$\vec{r}_{(7)} = (\frac{x^{\tilde{i}}_{(1)} + x^{\tilde{i}}_{(2)}}{2})\vec{e}_{\tilde{i}}; \qquad (11)$$

$$\vec{r}_{(9)} = (\frac{x^{\tilde{i}}_{(3)} + x^{\tilde{i}}_{(4)}}{2})\vec{e}_{\tilde{i}}$$

Eqns. (12), (13) describes how the natural values of the poles are scaled depending on the position of the poles $p_{(5)}$, $p_{(6)}$ in the positive or in the negative direction relating to the geometric center (g).

$$\overline{\theta}^1_{(5)} = |\vec{r}_{(5)} - \vec{r}_{(g)}|/|\vec{r}_{(8)} - \vec{r}_{(g)}| \text{ or } \theta^1_{(5)} = -|\vec{r}_{(5)} - \vec{r}_{(g)}|/|\vec{r}_{(10)} - \vec{r}_{(g)}|; \overline{\theta}^2_{(5)} = 0 \qquad (12)$$

$$\overline{\theta}^2_{(6)} = |\vec{r}_{(6)} - \vec{r}_{(g)}|/|\vec{r}_{(9)} - \vec{r}_{(g)}| \text{ or } \theta^1_{(5)} = -|\vec{r}_{(6)} - \vec{r}_{(g)}|/|\vec{r}_{(7)} - \vec{r}_{(g)}|; \overline{\theta}^1_{(6)} = 0 \qquad (13)$$

where $\vec{r}_{(g)}$ is the position vector of the geometric center. After computing the scaled values eqns. (12) and (13) representing the natural coordinates of the poles, respectively, the matrix in eq. (2) will be extended to include these values.

$$\theta^i_{(p)} = \begin{bmatrix} \theta^1_{(1)} & \theta^2_{(1)} \\ \theta^1_{(2)} & \theta^2_{(2)} \\ \theta^1_{(3)} & \theta^2_{(3)} \\ \theta^1_{(4)} & \theta^2_{(4)} \\ \theta^1_{(5)} & \theta^2_{(5)} \\ \theta^1_{(6)} & \theta^2_{(6)} \end{bmatrix} = \begin{bmatrix} -1 & -1 \\ +1 & -1 \\ +1 & +1 \\ -1 & +1 \\ \overline{\theta}^1_{(5)} & 0 \\ 0 & \overline{\theta}^2_{(6)} \end{bmatrix} \qquad (14)$$

Now, all the Cartesian and the natural values of the coordinates are prepared for the interpolation procedure and placing a complete two dimensional complete Pascal function of second order over the six nodal points. Substituting the special natural coordinates of the nodal points of eqn. (14) in eqn.(4) yields:

$$x^{\tilde{i}}_{(p)} - A^{(q)}_{(p)} a^{\tilde{i}}_{(q)} \qquad (15)$$

Where $A^{(q)}_{(p)}$ is a 6x6-matrix given by the following relation:

$$\begin{bmatrix} x_{(1)}^{\tilde{1}} & x_{(1)}^{\tilde{2}} \\ x_{(2)}^{\tilde{1}} & x_{(2)}^{\tilde{2}} \\ x_{(3)}^{\tilde{1}} & x_{(3)}^{\tilde{2}} \\ x_{(4)}^{\tilde{1}} & x_{(4)}^{\tilde{2}} \\ x_{(5)}^{\tilde{1}} & x_{(5)}^{\tilde{2}} \\ x_{(6)}^{\tilde{1}} & x_{(6)}^{\tilde{2}} \end{bmatrix} = \begin{bmatrix} 1 & -1 & -1 & 1 & 1 & 1 \\ 1 & 1 & -1 & 1 & -1 & 1 \\ 1 & 1 & 1 & 1 & 1 & 1 \\ 1 & -1 & 1 & 1 & -1 & 1 \\ 1 & \bar{\theta}_{(5)}^{1} & 0 & (\bar{\theta}_{(5)}^{1})^2 & 0 & 0 \\ 1 & 0 & \bar{\theta}_{(6)}^{2} & 0 & 0 & (\bar{\theta}_{(6)}^{2})^2 \end{bmatrix} \begin{bmatrix} a_{(1)}^{\tilde{1}} & a_{(1)}^{\tilde{2}} \\ a_{(2)}^{\tilde{1}} & a_{(2)}^{\tilde{2}} \\ a_{(3)}^{\tilde{1}} & a_{(3)}^{\tilde{2}} \\ a_{(4)}^{\tilde{1}} & a_{(4)}^{\tilde{2}} \\ a_{(5)}^{\tilde{1}} & a_{(5)}^{\tilde{2}} \\ a_{(6)}^{\tilde{1}} & a_{(6)}^{\tilde{2}} \end{bmatrix} \quad (16)$$

Solving the algebraic system of eqns. (16) for the undetermined parameters gives the geometric meaning of them:

$$a_{(q)}^{\tilde{i}} = B_{(q)}^{(p)} \, x_{(p)}^{\tilde{i}} \quad (17)$$

where $B_{(q)}^{(p)}$ is the inverse matrix of $A_{(p)}^{(q)}$. The existence of the solution requires a positive value of the determinant of the coefficients matrix ($\det(A_{(p)}^{(q)}) > 0$).

Substituting this result in eqns. (4) leads to the following relationship between the Cartesian coordinates of an arbitrary point of the domain and the Cartesian coordinates of the nodal points.

$$x^{\tilde{i}} = N^{(q)} \, a_{(q)}^{\tilde{i}} = N^{(q)} \, B_{(q)}^{(p)} \, x_{(p)}^{\tilde{i}} \quad (18)$$

Note that $N^{(q)} B_{(q)}^{(p)}$ are the so called shape functions and eq. (18) has the same form as eq.(4). The parametric set $a_{(p)}^{\tilde{i}}$ is defined by eq. (17) and is dependent only on the constant known values of the Cartesian coordinates. The transformation relation before the geometric interpolation has also the same form after the geometric interpolation, except that the undetermined parameters $a_{(p)}^{\tilde{i}}$ are replaced by the so called generalized parameters $B_{(q)}^{(p)} x_{(p)}^{\tilde{i}}$.

The relation (18) is in general a non linear system of two equations between the Cartesian coordinates $(x^{\tilde{1}}, x^{\tilde{2}})$ and the natural coordinates $(\theta^1, \theta^2)$. The Cartesian coordinates can be computed straightforward assuming known natural coordinates but if the Cartesian coordinates are given the system of equations become nonlinear owing to the existence of the non linear terms.

In the special cases, when two opposite edges, or every two of them are parallel, the solution of the algebraic system of eq. (16) includes the following results for the undetermined parameters $a_{(4)}^1, a_{(4)}^2 \, ; \, a_{(6)}^1, a_{(6)}^2$ associated with the quadratic natural variables $(\theta^1)^2 \, ; \, (\theta^2)^2$:

$$\bar{\theta}_{(5)}^1 \to \infty ; a_{(4)}^1 = 0, \, a_{(4)}^2 = 0 \quad (19)$$

$$\bar{\theta}_{(6)}^2 \to \infty; a_{(6)}^1 = 0, a_{(6)}^2 = 0 \tag{20}$$

$$\bar{\theta}_{(5)}^1 \to \infty, \bar{\theta}_{(5)}^1 \to \infty ; a_{(4)}^1 = 0, a_{(4)}^2 = 0 ; a_{(6)}^1 = 0, a_{(6)}^2 = 0 \tag{21}$$

For the last case (eq. (21)), the Pascal scheme delivers the Lagrangian scheme.

Now, the following geometric properties are defined, for more details see for example [23], [24], [25]:

- The position vector of an arbitrary point of the domain:

$$\vec{r} = x^{\tilde{i}} \vec{e}_{\tilde{i}} \tag{22}$$

- Covariant basis vectors

$$\vec{g}_\alpha = \vec{r}_{,\alpha} = x^{\tilde{i}}{}_{,\alpha} \vec{e}_{\tilde{i}} \tag{23}$$

- Derivatives of the covariant basis vectors

$$\vec{g}_{\alpha,\beta} = \vec{r}_{,\alpha\beta} = x^{\tilde{i}}{}_{,\alpha\beta} \vec{e}_{\tilde{i}} \tag{24}$$

- Metric coefficients

$$g_{\alpha\beta} = \vec{g}_\alpha \cdot \vec{g}_\beta = \vec{r}_{,\alpha} \cdot \vec{r}_{,\beta} = x^{\tilde{i}}{}_{,\alpha} \cdot x^{\tilde{j}}{}_{,\beta} \vec{e}_{\tilde{i}} \vec{e}_{\tilde{j}} \tag{25}$$

- Contra-variant metric coefficients

$$g^{\alpha\beta} = (g_{\alpha\beta})^{-1} \tag{26}$$

- Contra-variant basis vectors

$$\vec{g}^\alpha = g^{\alpha\beta} \cdot \vec{g}_\beta \tag{27}$$

- Chrisstoffel symbols of the first kind

$$\Gamma_{\alpha\beta\gamma} = g^{i}_{\alpha,\beta} \cdot g^{j}_{\gamma} \cdot \delta_{ij} \tag{28}$$

- Chrisstoffel symbols of the second kind

$$\Gamma^{\gamma}_{\alpha\beta} = g^{i}_{\alpha,\beta} \cdot g^{\gamma}_{i} \tag{29}$$

It is found that the following relations apply

$$\vec{g}_{\alpha,\beta} = \Gamma_{\alpha\beta\gamma} \cdot \vec{g}^{\gamma}$$
$$\vec{g}_{\alpha,\beta} = \Gamma^{\gamma}_{\alpha\beta} \cdot \vec{g}_{\gamma}$$
(30)

Note that the differential geometry properties at an arbitrary domain point changes significantly from that produced by a Lagrangian scheme owing to the existence of the additional quadratic natural variables $(\theta^1)^2$ ; $(\theta^2)^2$ in the presented Pascal scheme. The derivatives of the position vector of an arbitrary points includes also more terms than that produced by the Lagtangian scheme. It is expected that this change in the differential geometry properties must affect the results of numerical analysis positively. The use of coordinate system $(x^{*1}, x^{*2})$ located at the geometric centre or a rotated coordinate system $(x^1, x^2)$ defined from the directions of the covariant and contra variant basis vectors, see Fig.1, enhances the calculation of some differential geometry properties listed above.

## 3  2D THIRD ORDER INTERPOLATION SCHEME

In this case, a ten term Pascal Polynomial is used in order to approximate the Cartesian variables $x^{\tilde{1}}$ , $x^{\tilde{2}}$ as follows:

$$\begin{bmatrix} x^{\tilde{1}} & x^{\tilde{2}} \end{bmatrix} = \begin{bmatrix} 1 & \theta^1 & \theta^2 & (\theta^1)^2 & \theta^1\theta^2 & (\theta^2)^2 & (\theta^1)^3 & (\theta^1)^2\theta^2 & \theta^1(\theta^2)^2 & (\theta^2)^3 \end{bmatrix} \begin{bmatrix} a^{\tilde{1}}_{(1)} & a^{\tilde{2}}_{(1)} \\ a^{\tilde{1}}_{(2)} & a^{\tilde{2}}_{(2)} \\ a^{\tilde{1}}_{(3)} & a^{\tilde{2}}_{(3)} \\ a^{\tilde{1}}_{(4)} & a^{\tilde{2}}_{(4)} \\ a^{\tilde{1}}_{(5)} & a^{\tilde{2}}_{(5)} \\ a^{\tilde{1}}_{(6)} & a^{\tilde{2}}_{(6)} \\ a^{\tilde{1}}_{(7)} & a^{\tilde{2}}_{(7)} \\ a^{\tilde{1}}_{(8)} & a^{\tilde{2}}_{(8)} \\ a^{\tilde{1}}_{(9)} & a^{\tilde{2}}_{(9)} \\ a^{\tilde{1}}_{(10)} & a^{\tilde{2}}_{(10)} \end{bmatrix}$$
(31)

The four midpoints of the edges in addition to the six nodal point from the 2D second order approximation case are used as significant nodal points. The corresponding natural coordinates for these points are:

$$\begin{bmatrix} \theta^1_{(1)} & \theta^3_{(1)} \\ \theta^1_{(2)} & \theta^3_{(2)} \\ \theta^1_{(3)} & \theta^3_{(3)} \\ \theta^1_{(4)} & \theta^3_{(4)} \\ \theta^1_{(5)} & \theta^3_{(5)} \\ \theta^1_{(6)} & \theta^2_{(6)} \\ \theta^1_{(7)} & \theta^2_{(7)} \\ \theta^1_{(8)} & \theta^2_{(8)} \\ \theta^1_{(9)} & \theta^2_{(9)} \\ \theta^1_{(10)} & \theta^2_{(10)} \end{bmatrix} = \begin{bmatrix} -1 & -1 \\ +1 & -1 \\ +1 & +1 \\ -1 & +1 \\ \bar{\theta}^1_{(5)} & 0 \\ 0 & \bar{\theta}^2_{(6)} \\ 0 & -1 \\ +1 & 0 \\ 0 & +1 \\ -1 & 0 \end{bmatrix} \qquad (32)$$

$\bar{\theta}^1_{(5)}$ the scaled natural coordinate of the pole $p_{(5)}$, are now computed as a ratio between the length of the curves (g)(8) and (g)(5) or , (g)(10), (g)(5)..

$\bar{\theta}^2_{(6)}$ the scaled natural coordinates of the poles $p_{(6)}$ , are now computed as a ratio between the length of the curves (g)(9) and (g)(6) or , (g)(7), (g)(6).

The interpolation procedure leads to a 10x10 system of equations, which can be solved to find the undetermined parameters $_a\tilde{i}_{(p)}$ dependent on the Cartesian coordinates of the nodal points. It can be shown that the inverse matrix of $A^{(q)}_{(p)}$ exists even for straight edges, which means that the first interpolation scheme is included in last higher order case. The higher order scheme can be used for elements with curved edges represented by 1D two order polynomials.

Finally, a fifteen term Pascal Polynomial for approximating the geometry of quadrilateral element with edges represented by 1D third order polynomials is currently under investigation.

The fifteen free parameters are determined from the natural coordinates of the four nodal points (vertices) of the quadrilateral element and the two intersections points of the lines crossing every two opposite edges (poles) as well as of eight nodal points placed on the four edges and those of the geometric centre of the element

## 4 NUMERICAL COMPARISION

Example: A quadrilateral domain is defined by its four vertices [(1), (2) (3), (4)] (nodal points)

\

$$x^{\tilde{i}}_{(p)} = \begin{matrix} (p) \downarrow \\ \phantom{x} \end{matrix} \begin{bmatrix} 1. & 1. \\ 4. & 2. \\ 3. & 4. \\ 2. & 5. \end{bmatrix} \xrightarrow{\tilde{i}} \tag{33}$$

Solving the equations of every two opposite edges of the quadrilateral shape simultaneously delivers the Cartesian coordinates of its poles. Updating $x^{\tilde{i}}_{(p)}$ with the results gives the extended matrix of it

$$x^{\tilde{i}}_{(p)} = \begin{bmatrix} 1.0000 & 1.000 \\ 4.0000 & 2.000 \\ 3.000 & 4.000 \\ 2.000 & 5.000 \\ 4.7500 & 2.250 \\ 2.1667 & 5.6667 \end{bmatrix} \tag{34}$$

The corresponding natural coordinates of the nodal points are

$$\theta^{i}_{(p)} = \begin{bmatrix} \theta^1_{(1)} & \theta^2_{(1)} \\ \theta^1_{(2)} & \theta^2_{(2)} \\ \theta^1_{(3)} & \theta^2_{(3)} \\ \theta^1_{(4)} & \theta^2_{(4)} \\ \bar{\theta}^1_{(5)} & 0 \\ 0 & \bar{\theta}^2_{(6)} \end{bmatrix} = \begin{bmatrix} -1 & -1 \\ +1 & -1 \\ +1 & +1 \\ -1 & +1 \\ 3.80788655 & 0 \\ 0 & 1.76698110 \end{bmatrix} \tag{35}$$

$\bar{\theta}^1_{(5)}, \bar{\theta}^2_{(6)}$ are computed by applying eqns. (12) and (13).

Following the computational path from (15) to (18) yields the generalized coordinates:

$$a^{\tilde{i}}_{(p)} = \begin{bmatrix} 1.78769652 & 4.48213067 \\ 1.00000000 & 0.00000000 \\ 0.00000000 & 1.50000000 \\ -0.23647699 & 0.08169214 \\ -0.50000000 & -0.50000000 \\ 0.94878047 & -1.56382281 \end{bmatrix} \tag{36}$$

and the transformation relation between natural and Cartesian coordinates is then given by:

$$[x^{\tilde{1}} \quad x^{\tilde{2}}] = [1 \quad \theta^1 \quad \theta^2 \quad (\theta^1)^2 \quad \theta^1\theta^2 \quad (\theta^2)^2] \begin{bmatrix} 1.78769652 & 4.48213067 \\ 1.00000000 & 0.00000000 \\ 0.00000000 & 1.50000000 \\ -0.23647699 & 0.08169214 \\ -0.50000000 & -0.50000000 \\ 0.94878047 & -1.56382281 \end{bmatrix} \quad (37)$$

Substituting the natural coordinates of the four vertices and the two poles yields to the Cartesian coordinates of these points. Substituting also the natural equation of an edge gives information about how this edge is bounded. It can be easily proved that the essential properties of shape functions are preserved.

The Lagrangian scheme for this transformation is:

$$[x^{\tilde{1}} \quad x^{\tilde{2}}] = [1 \quad \theta^1 \quad \theta^2 \quad \theta^1\theta^2] \begin{bmatrix} 2.50 & 3.00 \\ 1.00 & 0.00 \\ 0.00 & 1.50 \\ -0.5 & -0.5 \end{bmatrix} \quad (38)$$

The covariant basis vectors for the Pascal scheme are now given by:

$$\vec{g}_\alpha = g_\alpha^{\tilde{i}} \vec{e}_{\tilde{i}} \quad ; \quad g_\alpha^{\tilde{i}} = \begin{bmatrix} g_1^{x^{\tilde{1}}} & g_1^{x^{\tilde{2}}} \\ g_2^{x^{\tilde{1}}} & g_2^{x^{\tilde{2}}} \end{bmatrix} = \begin{bmatrix} 0 & 1 & 0 & 2\theta^1 & \theta^2 & 0 \\ 0 & 0 & 1 & 0 & \theta^1 & 2\theta^2 \end{bmatrix} \begin{bmatrix} 1.78769652 & 4.48213067 \\ 1.00000000 & 0.00000000 \\ 0.00000000 & 1.50000000 \\ -0.23647699 & 0.08169214 \\ -0.50000000 & -0.50000000 \\ 0.94878047 & -1.56382281 \end{bmatrix} \quad (39)$$

Note that, the determinant of the Jacobian at the geometric centre is equal to quarter of the area of the quadrilateral element. The integration of the Jacobian determinant over the interval [-1,+1] gives also the area of the element.

The Lagrangian scheme gives:

$$\vec{g}_\alpha = g_\alpha^{\tilde{i}} \vec{e}_{\tilde{i}} \quad ; \quad g_\alpha^{\tilde{i}} = \begin{bmatrix} g_1^{x^{\tilde{1}}} & g_1^{x^{\tilde{2}}} \\ g_2^{x^{\tilde{1}}} & g_2^{x^{\tilde{2}}} \end{bmatrix} = \begin{bmatrix} 0 & 1 & 0 & \theta^2 \\ 0 & 0 & 1 & \theta^1 \end{bmatrix} \begin{bmatrix} 2.50 & 3.00 \\ 1.00 & 0.00 \\ 0.00 & 1.50 \\ -0.5 & -0.5 \end{bmatrix} \quad (40)$$

Note that, six generalized parameters are identical in both schemes. They are the four components of the basis vectors in the geometric centre and two components of their derivatives. At an arbitrary point of the quadrilateral element the functions defining these vectors are different.

$$\vec{g}_{\alpha,\beta} = g_\alpha^{\tilde{i}} e_{\tilde{i}} \quad ; \quad g_{\alpha,\beta}^{\tilde{i}} = \begin{array}{|c|c|} \hline g_{1,1}^{x^{\tilde{1}}} & g_{1,1}^{x^{\tilde{2}}} \\ \hline g_{2,1}^{x^{\tilde{1}}} & g_{2,1}^{x^{\tilde{2}}} \\ \hline g_{1,2}^{x^{\tilde{1}}} & g_{1,2}^{x^{\tilde{2}}} \\ \hline g_{2,2}^{x^{\tilde{1}}} & g_{2,2}^{x^{\tilde{2}}} \\ \hline \end{array} = \begin{bmatrix} 0 & 0 & 0 & 2 & 0 & 0 \\ 0 & 0 & 0 & 0 & 1 & 0 \\ 0 & 0 & 0 & 0 & 1 & 0 \\ 0 & 0 & 0 & 0 & 0 & 2 \end{bmatrix} \begin{bmatrix} 1.78769652 & 4.48213067 \\ 1.00000000 & 0.00000000 \\ 0.00000000 & 1.50000000 \\ -0.23647699 & 0.08169214 \\ -0.50000000 & -0.50000000 \\ 0.94878047 & -1.56382281 \end{bmatrix} \quad (41)$$

All components of the derivatives of the covariant basis vectors produced by Pascal scheme are presented, whereas there are four components of them equal to zero by applying the Lagrangian scheme.

$$\vec{g}_{\alpha,\beta} = g_\alpha^{\tilde{i}} e_{*_{\tilde{i}}} \quad ; \quad g_{\alpha,\beta}^{\tilde{i}} = \begin{array}{|c|c|} \hline g_{1,1}^{x^{\tilde{1}}} & g_{1,1}^{x^{\tilde{2}}} \\ \hline g_{2,1}^{x^{\tilde{1}}} & g_{2,1}^{x^{\tilde{2}}} \\ \hline g_{1,2}^{x^{\tilde{1}}} & g_{1,2}^{x^{\tilde{2}}} \\ \hline g_{2,2}^{x^{\tilde{1}}} & g_{2,2}^{x^{\tilde{2}}} \\ \hline \end{array} = \begin{bmatrix} 0 & 0 & 0 & 0 \\ 0 & 0 & 0 & 1 \\ 0 & 0 & 0 & 1 \\ 0 & 0 & 0 & 0 \end{bmatrix} \begin{bmatrix} 2.50 & 3.00 \\ 1.00 & 0.00 \\ 0.00 & 1.50 \\ -0.5 & -0.5 \end{bmatrix} \quad (42)$$

## 5  CONCLUSIONS

A complete Pascal interpolation scheme for approximating the differential geometry properties of a quadrilateral shape and establishing a transformation relation between natural and Cartesian coordinates has been presented. It is found that the presented scheme is more general and approximate the geometry of the quadrilateral shape completely. It recovers the Lagrangian scheme in the special case, when every two opposite edges of the shape are parallel. The Lagrangian scheme for approximating the geometry of rectangular and rhombic shapes appears to be a special case of the presented approach.

To Cite this article,

[]   Abo Diab, S., Complete Pascal Interpolation Scheme For Approximating The Geometry Of A Quadrilateral Element, *Department of Structural Mechanics, Tishreen University, Faculty of Civil Engineering, Lattakia, Syria( 2017)*.